\documentclass[12pt,a4paper]{amsart}

\newcommand{\N}{\rm{I\!N}}         
\newcommand{\R}{\rm{I\!R}}

\newtheorem{theorem}{Theorem}
\newtheorem{definition}[theorem]{Definition}

\newtheorem{lemma}[theorem]{Lemma}
\newcommand{\eps}{\varepsilon}

\newcommand{\norm}[1]{\|#1\|}

\newcommand{\betr}[1]{| #1  |}

\newcommand{\ebew}{\hfill{\rule{1.2ex}{1.2ex}}}
\newcommand{\bgl}{\begin{eqnarray}}
\newcommand{\bglst}{\begin{eqnarray*}}
\newcommand{\egl}{\end{eqnarray}}
\newcommand{\eglst}{\end{eqnarray*}}

\newcommand{\Ref}[1]{(\ref{#1})}

\newcommand{\mdE}{\;|}
\newcommand{\ball}{{\mathcal B}}
\newcommand{\be}{\begin{enumerate}}
\newcommand{\ee}{\end{enumerate}}
\renewcommand{\limsup}{\overline{\lim}\;}
\renewcommand{\liminf}{\underline{\lim}\;}

\begin{document}
\title{Boundaries for Banach spaces determine weak compactness}
\author{Hermann Pfitzner}
\address{Universit\'e d'Orl\'eans\\
BP 6759\\
F-45067 Orl\'eans Cedex 2\\France}
\email{hermann.pfitzner@univ-orleans.fr}
\begin{abstract}
A boundary for a Banach space is a subset of the dual unit sphere with the property that each element of the Banach space attains its norm on
an element of that subset. Trivially, the pointwise convergence with respect to such a boundary is coarser than the weak topology on the Banach space.
Godefroy's Boundary Problem asks whether nevertheless both topologies have the same bounded compact sets.
This paper contains the answer in the positive.
\end{abstract}
\maketitle\noindent
This paper deals with boundaries for Banach spaces in the sense of
\begin{definition}
Let $X$ be a real Banach space. A subset $B$ of ${\mathcal S}(X^*)$ is called a boundary for $X$ if for each $x\in X$ there
is $b\in B$ such that $b(x)=\norm{x}$.
\end{definition}
\noindent
The set of extremal points of the dual unit ball of a Banach space is a well-known example of a boundary.
In 1980 Bourgain and Talagrand \cite{BourgTal} showed that a norm bounded subset of a Banach space is weakly compact if it is compact in the pointwise topology on the
set of extremal points of the dual unit ball.
Some years later Godefroy  \cite{G-boundaries}  asked whether the result still holds if the extremal points are replaced by an arbitrary boundary.

Theorem \ref{Satz6} shows that the answer is ``yes''.
The proof can be described vaguely as an
amalgam of 
Behrends' quantitative version of Rosenthal's $l^1$-Theorem (cf.\ Lemma \ref{Lemma1}),
Simons' equality (cf.\ Lemma \ref{LemmaSimons}),
a variant of Hagler-Johnson's construction (cf.\ Lemma \ref{Hauptlemma})
and James' distortion theorem (everywhere).

Besides the result of Bourgain and Talagrand a positive answer to Godefroy's question has been known in the important case when the set in question is convex \cite[p.\ 44]{G-Murcia}.
Bourgain's and Talagrand's proof relies on  the results of \cite{BourgFremTal} which are of topological nature and do not seem applicable here
because general boundaries seem to lack sufficient topological structures. As a - technically quite different - substitute we use Simons' equality \cite{Sim-Fatou}
(or, if the reader prefers, Simons' inequality \cite{Sim-boundary}, see \cite[Th. 3.48]{FHHMPZ})
which has been advocated at several instances by Godefroy (e. g. \cite{Cas-Go, Go-SimIneq}).
The idea to look for the key Lemma \ref{Hauptlemma} of the present proof was inspired by the result of Cascales et al.  \cite{Cas-Man-Ver, Cas-Shv} on the existence of
an independent sequence.
\smallskip\\
Throughout this article $X$ denotes a  real Banach space, $X^*$ its dual, $\ball(X)$ its unit ball and ${\mathcal S}(X)$ its unit sphere.
The norm closed linear span of a subset $A$ of $X$ is written $[A]$. $\N$ starts at 1.
Our references for unexplained Banach space notions are \cite{JohLin} (for boundaries in particular see Ch.\ 15, Infinite Dimensional Convexity by Fonf, Lindenstrauss and Phelps),
\cite{Die-Seq} and \cite{LiTz12}.
\medskip\\
{\sc Acknowledgement}
It is a pleasure to thank G. Godefroy for having introduced me to the subject and for encouraging conversations.
\bigskip\\
To each bounded sequence $(x_{n})$ in $X$ we associate its James constant
\bglst
\eps_J (x_{n}) =\sup_m \inf_{\sum_{n\geq m}\betr{\alpha_n} =1}\norm{\sum_{n\geq m} \alpha_n x_n}.
\eglst
If a sequence is equivalent to the canonical basis of $l^1$ we call it simply an $l^1$-sequence.
Clearly, $\eps_J\geq0$ with $\eps_J (x_{n})>0$ if and only if there is an integer $m$ such that $(x_{n})_{n\geq m}$ is an $l^1$-sequence.

Two more moduli will be of importance. Let $D$ be a subset of $X^*$. We define
\bglst
\delta_D (x_{n}) = \sup_{x^*\in D}\big(\limsup x^* (x_{n}) -\liminf x^*(x_{n})\big),\quad
\delta_{\mathrm{HJ}, \; D} (x_{n}) = \sup_{x^*\in D}x^* (x_n).
\eglst
for bounded sequences $(x_n)$ and write $\delta =\delta_{\ball(X^*)}$ and $\delta_{\mathrm{HJ}} =\delta_{\mathrm{HJ}, \; \ball(X^*)}$ for short.
Clearly $\delta_D\geq0$ for all $D$ and $\delta_{\mathrm{HJ}}\geq 0$ and equality $\delta (x_{n})=0$ (respectively $\delta_{\mathrm{HJ}} (x_{n})=0$) holds if and only if
$(x_{n})$ is weakly Cauchy (respectively converges weakly to zero).
Both $\delta_D$ and $\delta_{\mathrm{HJ}, \; D}$ are non-increasing and similarly $\eps_J$ is non-decreasing if one passes to subsequences.
We introduce the notation $$\tilde{\delta}_D(x_n)=\inf_{n_k}\delta_D (x_{n_k}) \quad\mbox{ and }\quad
\tilde{\eps}_J(x_n)=\sup_{n_k}\eps_J (x_{n_k})$$
and say that $(x_{n})$ is $\delta_D$-stable if $\tilde{\delta}_D (x_{n})=\delta_D (x_{n})$. Likewise $(x_n)$ is $\eps_J$-stable
if $\tilde{\eps}_J(x_n)=\eps_J (x_n)$.

If one takes a norm preserving Hahn-Banach extension of the functional defined on $[x_n]_{n\geq m}$ by $x_n\mapsto (-1)^n\eps_J((x_k)_{k\geq m}) $ for all
$n\geq m$ and $m$ big enough then it is completely elementary but important to deduce that $\delta_{\mathrm{HJ}}  \geq \eps_J $
and $\delta \geq 2\eps_J$, even $\tilde{\delta} \geq 2\eps_J$.
While in general strict inequality may occur for $\delta_{\mathrm{HJ}}$
this cannot happen in $\tilde{\delta} \geq 2\eps_J$ as soon as $(x_n)$ is $\eps_J$-stable.
This follows from Behrends' quantitative version of Rosenthal's $l^1$-Theorem.
\begin{lemma}\label{Lemma1}
Let $X$ be a Banach space.
\be
\item\label{item1} Each bounded sequence in $X$ contains an $\eps_J$-stable subsequence.
\item\label{item2} Each bounded sequence in $X$ contains a $\delta$-stable subsequence.
\item\label{item3} If $(x_{n})$ is $\eps_J$-stable then
\bgl
\tilde{\delta}(x_n)=2\tilde{\eps}_J(x_n).\label{gla}
\egl
If moreover $(x_n)$ is $\delta$-stable then
\bgl
\tilde{\delta}(z_n)=2\tilde{\eps}_J(x_n).\label{glb}
\egl
whenever
\bglst
z_n=\frac{1}{2}\sum_{l=1}^{q}\lambda_l (x_{p_{2l}(n)}-x_{p_{2l-1}(n)})
\eglst
where $q\in\N$, $\sum_{l=1}^q \lambda_l=1$, $\lambda_l\geq0$ and where $(\{p_1(n), \ldots, p_{2q}(n)\})_{n\in\N}$ is a sequence of pairwise disjoint subsets
(of cardinality $2q$) of $\N$.
\ee
\end{lemma}\noindent
{\em Proof}:
\Ref{item1} The proof consists in a routine diagonal argument.
Choose a subsequence $(y_n^{(1)})$ of $(x_n)$ such that $\eps_J(y_n^{(1)})\geq \tilde{\eps}_J(x_n)-2^{-1}$.
Successiveley, choose a subsequence $(y_n^{(k+1)})$ of $(y_n^{(k)})$ for each $k\in\N$ such that
$\eps_J(y_n^{(k+1)})\geq \tilde{\eps}_J(y_n^{(k)})-2^{-(k+1)}$. Let $y_n=y_n^{(n)}$ and let $(z_n)$ be
a subsequence of $(y_n)$. Then
$\eps_J(y_n)\geq \eps_J(y_n^{(k+1)})\geq \tilde{\eps}_J(y_n^{(k)})-2^{-(k+1)}\geq \eps_J(z_n)-2^{-(k+1)}$
whence $\eps_J(y_n)\geq \eps_J(z_n)$ and $\eps_J(y_n)\geq \tilde{\eps}_J(y_n)$ that is $\eps_J(y_n)=\tilde{\eps}_J(y_n)$.
The proof of \Ref{item2} works alike.\\
\Ref{item3} Before the statement of the lemma we have already noticed that $\tilde{\delta}\geq 2\eps_J$ whence $\tilde{\delta}\geq 2\tilde{\eps}_J$
because $(x_n)$ is assumed to be $\eps_J$-stable. In order to prove the other inequality of \Ref{gla}
we may suppose that $\tilde{\delta}(x_n)>0$ because otherwise $\tilde{\delta}(x_n)=0$ and \Ref{gla} holds trivially.
In our notation Behrends' main result \cite[Th.\ 3.2]{Beh-ros} reads: If $\tilde{\delta}(x_n)>0$ then, given $\eta>0$,
there is a subsequence whose $\eps_J$-value is greater than $-\eta+(\tilde{\delta}(x_n))/2$.
Now let $\eta$ run through a sequence tending to zero and pass successively to according subsequences;
the process results in a diagonal sequence $(x_{n_k})$ such that $\eps_J(x_{n_k})\geq(\tilde{\delta}(x_n))/2$.
Hence $\tilde{\eps}_J(x_n)\geq(\tilde{\delta}(x_n))/2$ which proves \Ref{gla}.

For ``$\geq$'' of \Ref{glb} proceed as before the statement of the lemma
by using functionals of the kind $x_{p_{2l-i}(n)}\mapsto (-1)^i \eps_J((x_k)_{k\geq m})$, $i\in\{0,1\}$, $1\leq l\leq q$, with $n$ and $m$ big enough.
For the other inequality of \Ref{glb} note that $\delta$-stability and \Ref{gla} reduce it to $\tilde{\delta}(z_n)\leq \delta(x_n)$ which in turn by
subadditivity of $\delta$ is reduced to $\delta(x_{n_k}-x_{m_k})\leq 2\delta(x_n)$ where $(x_{n_k})$ and $(x_{m_k})$ are two disjoint subsequences of $(x_n)$.
But the latter inequality follows from the fact that for each $x^*\in \ball(X^*)$ both $\limsup x^*(x_{n_k}-x_{m_k})$ and $-\liminf x^*(x_{n_k}-x_{m_k})$
are each the difference of two cluster points of $(x^*(x_n))$ and hence majorized by $\delta(x_n)$.
\ebew\bigskip

If $B$ is a boundary for $X$, topological notions that refer to the $\sigma(X,B)$-topology are preceded by ``$B$-``, for example $B$-closed, $B$-compact.\\
The following lemma is a consequence of Simons' equality \cite{Sim-Fatou} which in our notation reads $\delta_{\mathrm{HJ}, B}=\delta_{\mathrm{HJ}}$.
\begin{lemma}\label{LemmaSimons}
Let $B$ be a boundary for $X$.
Then
\bgl
\delta_B=\delta \;\;\mbox{ and }\;\;\tilde{\delta}_B=\tilde{\delta}.\label{gl5a}
\egl
Moreover, if $(x_n)$ is an $\eps_J$-stable bounded sequence and $z_n$ is as in Lemma \ref{Lemma1} then
\bgl
\tilde{\delta}_B(z_n)=2\tilde{\eps}(x_n).\label{gl5b}
\egl
\end{lemma}\noindent
{\em Proof}: \Ref{gl5b} follows from \Ref{glb} and \Ref{gl5a} and the second half of \Ref{gl5a} is immediate from the first one which,
in turn, is a routine consequence of \Ref{glb} and Simons' equality: Fix $x^*\in\ball(X^*)$
and choose subsequences $(u_k)=(x_{n_k})$ and $(v_k)=(x_{m_k})$ such that
\bglst
\lefteqn{\limsup x^*(x_{n}) - \liminf x^*(x_{n})}\\
&=&\lim x^*(u_k-v_k)
                    \leq \delta_{\mathrm{HJ}}(u_k-v_k)=\delta_{\mathrm{HJ}, B}(u_k-v_k)\\
&\leq& \sup_{b\in B}\big(\limsup b(u_k)-\liminf b(v_k)\big)
     \leq \sup_{b\in B}\big(\limsup b(x_{n}) - \liminf b(x_{n})\big).
\eglst
This shows ''$\geq$`` of $\delta_B=\delta$ whereas ''$\leq$`` is trivial.\ebew\bigskip\\
Let $S=\bigcup_n S_n$ where $S_n=\{0, 1\}^n$ for each $n\in\N$.
If $\sigma\in S_n$ we write $\sigma=(\sigma_1, \ldots, \sigma_n)$ (with, of course, $\sigma_k\in\{0, 1\}$)
and for $i\in\{0, 1\}$ we write $(\sigma, i)$ for the element $(\sigma_1, \ldots, \sigma_n, i)\in S_{n+1}$.
Recall that a tree of non empty subsets of $\N$ is a sequence $(\Omega_{\sigma})_{\sigma\in S}$ such that
$\Omega_{(\sigma,0)}$ and $\Omega_{(\sigma,1)}$ are two disjoint non empty (hence infinite) subsets of $\Omega_{\sigma}\subset \N$ for all $\sigma\in S$.
\begin{lemma}\label{Hauptlemma}
Let $B$ be a boundary for $X$, $(x_n)$ be an $l^1$-sequence and $\eta_k>0$ decrease to zero.
Then there are a sequence $(b_k)$ in $B$, a tree $(\Omega_{\sigma})_{\sigma\in S}$ and $\eps\geq \eps_J(x_n)$ such that for all $k$
\bgl
b_k(x_n-x_{n'})>2(1-\eta_k)\eps \;\;\mbox{ if }\; n\in \Omega_{\sigma}, n'\in \Omega_{\sigma'}, \;\mbox{}\;\sigma, \sigma'\in S_k\mbox{ and } \sigma_k=0, \sigma'_k=1. \label{gl1}
\egl
Furthermore, if the set $\{x_n\mdE n\in\N\}$ is relatively $B$-compact in $X$, there is a sequence $(y_m)$ of $B$-cluster points of the $x_n$ such that
\bgl
b_k(y_m - y_{m'}) \geq 2(1-\eta_k)\eps\quad\mbox{ if }\quad m\leq k<m',\;\;k,m,m'\in\N.   \label{gl11}
\egl
\end{lemma}\noindent
{\em Proof}:
Choose $\Omega_{\emptyset}\subset\N$ such that the $x_n$ with indices in $\Omega_{\emptyset}$ form an $\eps_J$-stable and $\delta$-stable
$l^1$-sequence which exists by Lemma \ref{Lemma1}
and set $\eps=\eps_J((x_n)_{n\in\Omega_{\emptyset}})$. Then $\eps\geq \eps_J(x_n)$.\\
The $b_k$ and $\Omega_{\sigma}$ for $\sigma\in S_k$ will be constructed by induction over $k\in\N$.

For the first induction step $k=1$, equalities \Ref{gl5a} and \Ref{gla} allow to find $b_1\in B$ and a subsequence $(n_l)$ of $\Omega_{\emptyset}$ such that
$$b_1(x_{n_{2l}}-x_{n_{2l'+1}})>2(1-\eta_1)\eps$$
for all $l,l' \in\N$. It remains to set $\Omega_{(i)}=\{n_{2l+i}\mdE l\in\N\}$ for $i=0$ and $i=1$ in order to settle the first induction step.

Suppose that $b_1, \ldots, b_k$ and $\Omega_{\sigma}=\{\omega_{\sigma}(n)\mdE n\in\N\}$ for $\sigma\in S_k$ have been constructed
according to \Ref{gl1}. (Of course, we suppose the $\omega_{\sigma}:\N\rightarrow\N$ to be strictly increasing.)
Set $\eta=\eta_{k+1}/2^{k+1}$.
Apply \Ref{gl5b} to
$$z_{n}=2^{-k}\sum_{\sigma\in S_k} (-1)^{\sigma_k} \; x_{\omega_{\sigma}(n)}
$$
(with $q=2^{k-1}$ and $\lambda_l=1/q$ for $l\leq q$)
in order to get $b_{k+1}\in B$ and a sequence $(n_l)$ of integers such that
\bgl
b_{k+1}(z_{n_{2l}}-z_{n_{2l+1}})>2(1-\eta)\eps \quad\mbox{ for all }\quad l\in \N.\label{gl12}
\egl
Note that (by omitting at most finitely many members of the $n_k$) one has furthermore that
\bgl
b_{k+1}(x_{\omega_{\sigma}(n_{2l})} - x_{\omega_{\sigma'}(n_{2l+1})})<2(1+\eta)\eps\quad\mbox{ for all }\quad l\in \N, \sigma, \sigma'\in S_k\label{gl6}
\egl
because by \Ref{gla} the difference of two cluster points of $b_{k+1}(x_n)$ cannot exceed $\delta(x_n)=\tilde{\delta}(x_n)\leq 2\eps$.\\
Define, for all $\sigma\in S_k$,
\bglst
\left.\begin{array}{rl}
\omega_{(\sigma, 0)}(l)&=\omega_{\sigma}(n_{2l})\\
\omega_{(\sigma, 1)}(l)&=\omega_{\sigma}(n_{2l+1})
\end{array}\right\}    \mbox{if $\sigma_k=0$\quad and }
\left.\begin{array}{rl}
\omega_{(\sigma, 0)}(l)&=\omega_{\sigma}(n_{2l+1})\\
\omega_{(\sigma, 1)}(l)&=\omega_{\sigma}(n_{2l})
\end{array}\right\} \mbox{if $\sigma_k=1$.}
\eglst
With this notation we have
\bglst
z_{n_{2l}}-z_{n_{2l+1}}&=&2^{-k}\sum_{\sigma\in S_k}(-1)^{\sigma_k} (x_{\omega_{\sigma}(n_{2l})}-x_{\omega_{\sigma}(n_{2l+1})})\\
&=&2^{-k}\sum_{\sigma\in S_k}(x_{\omega_{(\sigma, 0)}(l)}-x_{\omega_{(\sigma, 1)}(l)}).
\eglst
Consider $\sigma, \sigma'\in S_k$ and distinguish the two cases $\sigma=\sigma'$ and $\sigma\neq\sigma'$.
In the first case we have
$$x_{\omega_{(\sigma, 0)}(l)}-x_{\omega_{(\sigma, 1)}(l)} =2^k (z_{n_{2l}}-z_{n_{2l+1}})-\sum_{\tau\in S, \tau\neq\sigma}(x_{\omega_{(\tau,0)}(l)}-x_{\omega_{(\tau,1)}(l)}),$$
and in the second case
$$x_{\omega_{(\sigma, 0)}(l)}-x_{\omega_{(\sigma',1)}(l)} =2^k (z_{n_{2l}}-z_{n_{2l+1}})-\big[(x_{\omega_{(\sigma',0)}(l)}-x_{\omega_{(\sigma,1)}(l)})+{\sum}' (x_{\omega_{(\tau,0)}(l)}-x_{\omega_{(\tau,1)}(l)})\big]$$
where the sum ${\sum}'$ runs over all $\tau\in S_k$ such that $\tau\neq\sigma$ and $\tau\neq\sigma'$.
In both cases we obtain
$$b_{k+1}(x_{\omega_{(\sigma, 0)}(l)}-x_{\omega_{(\sigma',1)}(l)}) > 2^{k+1} \, (1-\eta)\eps - 2\,(2^k -1)\,(1+\eta)\eps =2(1-\eta_{k+1}+\eta)\eps$$
for all $l\in\N$ by \Ref{gl12} and \Ref{gl6}.
Finally, since all $b_{k+1}(x_n)$ are contained in a compact subset of $\R$ there is an infinite set $N\subset\N$ such that
$$b_{k+1}(x_{\omega_{(\sigma, 0)}(l)}-x_{\omega_{(\sigma',1)}(l')}) > 2(1-\eta_{k+1})\eps$$
for all $l, l'\in N$.
It remains to set $\Omega_{(\sigma,i)}= \omega_{(\sigma, i)}(N)$ for $i\in\{0,1\}$ and  all $\sigma\in S_k$.
This shows \Ref{gl1} for $k+1$ and ends the induction.

For the last part of the lemma fix $m$, take, for all $k\geq m$,
\bglst
\sigma^{(k)}=(\underbrace{1,\ldots,1}_{m-1},\underbrace{0,\ldots,0}_{k-m+1})\in S_k,
\eglst
take $n_k$ to be the $k$-th element of $\Omega_{\sigma^{(k)}}$ and define $y_m$ to be a $B$-cluster point of the $x_{n_k}$.
Then, whenever $1\leq m\leq k<m'$, there are $\sigma, \sigma' \in S_k$ with $\sigma_k=0$, $\sigma'_k=1$
such that $y_m$(respectivlely $y_{m'}$) is a $B$-cluster point of the $x_n$
with indices in $\Omega_{\sigma}$ (respectively in $\Omega_{\sigma'}$)) and \Ref{gl11} follows from \Ref{gl1}.\ebew\bigskip

\begin{theorem}\label{Thl1}In a real Banach space with a boundary $B$, a bounded relatively $B$-compact set cannot contain an $l^1$-sequence.
\end{theorem}\noindent
{\em Proof}:
Suppose  to the contrary that there is a relatively $B$-compact $l^1$-sequence in a real Banach space $X$.
By Lemma \ref{Lemma1}  this $l^1$-sequence contains a $\delta$-stable subsequence which we denote by $(x_n)$.
Let $(\eta_k)$ be a decreasing sequence of
positive numbers with limit zero.
Take $\eps\geq\eps_J(x_n)>0$ and two sequences $(b_k)$  and $(y_m)$ that fulfill \Ref{gl11} of Lemma \ref{Hauptlemma} and set
$$x=\big(\sum 2^{-m}y_m\big)-y=\sum2^{-m}(y_m-y)$$
where $y$ is a $B$-cluster point of the $y_m$.
Inequality \Ref{gl11} extends to $y$ in the sense that $b_k(y_m-y) \geq 2(1-\eta_k)\eps$ for all $m\leq k$.
Note that $y$ is also a $B$-cluster point of the $x_n$.
Therefore the difference $b(y_m)-b(y)$ is majorized by $\delta(x_n)=\tilde{\delta}(x_n)\leq 2\eps$ for all $b\in B$ hence $\norm{x}\leq2\eps$.
But actually $\norm{x}=2\eps$: Let $\eta>0$, let $m_0$ such that $\norm{\sum_{m_0+1}^{\infty}2^{-m}(y_m-y)}<\eta$
and let $k>m_0$.
Then
\bglst
b_k(x)&\geq&\big(\sum_{m=1}^{m_0}2^{-m}b_k(y_m-y)\big)-\eta\\
&\geq& \big(\sum_{m=1}^{m_0}2^{-m}\;2(1-\eta_k)\eps)\big)-\eta\\
&=&2(1-\eta_k)(1-2^{-m_0})\eps-\eta
\eglst
which shows that $\sup_k b_k(x)\geq2\eps$ and thus
$$\norm{x}=\sup_k b_k(x)=2\eps.$$
Since $B$ is a boundary it contains a $b_0$ such that $b_0(x)=2\eps$.
So $b_0(y_m-y)=2\eps$ for all $m$ and $b_0(y_m)=2\eps+b_0(y)$. But $y$ is a $B$-cluster point of the $y_m$ thus
$b_0(y)=2\eps+b_0(y)$, a contradiction which ends the proof.\ebew\bigskip
\begin{theorem}\label{Satz6}
Let $B$ be a boundary for the real Banach space $X$. Then a $B$-compact bounded subset of $X$ is weakly compact.
\end{theorem}\noindent
{\em Proof}:
Let $A$ be a bounded $B$-compact subset of $X$.
By the theorem of Eberlein-\v{S}mulyan, in order to show that $A$ is weakly compact it is enough to prove that
each sequence in $A$ admits a weakly convergent subsequence.
Take a bounded sequence in $A$ and denote by $(x_n)$ the $\eps_J$-stable subsequence of it which exists by Lemma \ref{Lemma1}.
Theorem \ref{Thl1} entails that $\tilde{\eps}_J(x_n)=0$ hence $\tilde{\delta}(x_n)=0$ by \ref{gla}. That is, $(x_n)$ is weakly Cauchy hence $B$-Cauchy.
Since moreover $A$ is $B$-compact, $(x_n)$ $B$-converges (as a sequence) to a limit, say $x\in X$.
Now a final application of Simons' equality ends the proof because by
\bglst
\sup_{x^*\in\ball(X^*)}\limsup x^*(x_n -x) = \sup_{b\in B}\limsup b(x_n -x)=0
\eglst
we see that $(x_n)$ converges weakly to $x$.
\ebew
\end{document}